\documentclass[12pt]{amsart}
\usepackage{amssymb}

 \newtheorem{theorem}{Theorem}[section]
 \newtheorem{corollary}[theorem]{Corollary}
 \newtheorem{lemma}[theorem]{Lemma}
 \newtheorem{prop}[theorem]{Proposition}
 \theoremstyle{definition}
 
 \theoremstyle{remark}
 \newtheorem{remark}[theorem]{Remark}
 
 \numberwithin{equation}{section}

\newcommand{\be}{\begin{equation}}
\newcommand{\ee}{\end{equation}}

\def\Z{{\mathbb Z}}

\def\Tn{{\mathbb T}^n}
\def\Zn{{\mathbb Z}^n}

\def\supp{\operatorname{supp}}

\def\Rn{{{\mathbb R}^n}}
\def\Scal{\mathcal{S}}
\def\Fcal{\mathcal{F}}
\def\Dcal{\mathcal{D}}
\def\Lcal{\mathcal{L}}
\def\Ecal{\mathcal{E}}

\begin{document}

\title[On the Fourier analysis of operators on the torus]
{On the Fourier analysis of operators on the torus}

\author{Michael Ruzhansky}
\address{Michael Ruzhansky:
  Department of Mathematics \\
  Imperial College London \\
  180 Queen's Gate, London SW7 2AZ, United Kingdom}
\email{m.ruzhansky@imperial.ac.uk}

\thanks{The paper is the extended version of the talk given
by the authors at the $5^{th}$ ISAAC Congress, Catania, 2004.
\\ This first author would like to thank 
the UK Royal Society for its support.
The second author thanks the Academy of Finland
and Magnus Ehrnrooth Foundation for their support.}

\author{Ville Turunen} 
\address{Ville Turunen:
  Helsinki University of Technology\\ Institute of Mathematics\\
  P.O. Box 1100\\ FIN-02015 Finland\\ Finland}
\email{ville.turunen@hut.fi}

\subjclass{Primary 35L40; Secondary 58J40}

\keywords{Fourier integral operators, torus, pseudodifferential operators}

\begin{abstract}
Basic properties of Fourier integral operators on the torus
$\Bbb T^n=(\Bbb R/2\pi\Bbb Z)^n$ are studied by using
the global representations by Fourier series
instead of local representations.
The results can be applied to weakly
hyperbolic partial differential equations.
\end{abstract}

\maketitle

\section{Introduction}

In this paper we will discuss the version of the Fourier
analysis and pseudo-differential operators on the torus.
Using the toroidal Fourier transform we will show several
simplifications of the standard theory. We will also
discuss the corresponding toroidal version of Fourier
integral operators. To distinguish them from those defined
using the Euclidean Fourier transform, we will call them
{\em Fourier series operators}. The use of discrete
Fourier transform will allow to use global representation
of these operators, thus eliminating a number of topological
obstructions known in the standard theory. We will
prepare the machinery and describe how it can be further
used in the calculus of Fourier series operators and
applications to hyperbolic partial differential equations.
In fact, the form of the required discrete calculus 
is not a-priori clear, for example, the form of the discrete
Taylor's theorem best adopted to the calculus.
We will develop the corresponding version of the periodic
analysis similar in formulations to the standard Euclidean
theory.

It was realised already in the 1970s that on the torus,
one can study pseudodifferential operators globally
using Fourier series expansions,
in analogy to Euclidean pseudodifferential calculus.
These {\it periodic pseudodifferential operators}
were treated e.g. by Agranovich \cite{Agranovich1,Agranovich2}.
Contributions have been made by many authors, and
the following is a non-comprehensive list of the research on the torus:
Agranovich, crediting the idea to Volevich, proposed
the Fourier series representation of pseudodifferential operators.
Later, he proved the equivalence of the Fourier series representation
and H\"ormander's definition for $(1,0)$-symbol classes;
the case of classical pseudodifferential operators on the circle
had been treated by Saranen and Wendland \cite{SaranenWendland};
McLean \cite{McLean} proved the equivalence of
these approaches for H\"ormander's general $(\rho,\delta)$-classes,
by using charts;
in \cite{Turunen}, the equivalence for the case of $(1,0)$-classes
is proven by studying iterated commutators of pseudodifferential operators
and smooth vector fields.
Elschner \cite{Elschner} and Amosov \cite{Amosov}
constructed asymptotic expansions
for classical pseudodifferential operators;
these results were generalized for $(\rho,\delta)$-classes in
\cite{TurunenVainikko}.
There are plenty of papers considering applications
and numerical computation of pseudodifferential equations on torus,
e.g. spline approximations by
Pr\"ossdorf and Schneider \cite{ProssdorfSchneider},
physical applications by e.g.
Vainikko and Lifanov \cite{VainikkoLifanov1,
VainikkoLifanov2}, and many others.

On the other hand, the use of operators which are discrete in
the frequency variable allows one to weaken regularity 
assumptions on symbols with respect to $\xi$.
Symbols with low regularity in $x$ have been under intensive
study for many years, e.g.
Kumano-go and Nagase \cite{Kumanogonagase}, Sugimoto 
\cite{Sugimoto}, Boulkhemair \cite{Boulkhemair}, 
Garello and Morando \cite{GarelloMorando}, and many others.
However, in these papers one assumes symbols to be smooth
or sufficiently regular in $\xi$.
The discrete approach in this paper will allow us to reduce
regularity assumptions with respect to $\xi$. For example,
no regularity with respect to $\xi$ is assumed for $L^2$
estimates, and for elements of the calculus.
Moreover, one can consider scalar hyperbolic equations
with $C^1$ symbols with respect to $\xi$. For example,
this allows to construct parametrices to certain hyperbolic
systems with variable multiplicities. Details of such
constructions will appear in our forthcoming paper
\cite{RT}.

Let us now fix the notation.
Let ${\Scal}(\Bbb R^n)$ be the Schwartz test function space
with its usual topology,
and let ${\Scal}'(\Bbb R^n)$ be its dual,
the space of tempered distributions.
Let ${\Fcal}_E:{\Scal}(\Bbb R^n)\to{\Scal}'(\Bbb R^n)$
be the Euclidean Fourier transform
(hence the subscript $_E$)
defined by
$$
  ({\Fcal}_E f)(\xi) = \widehat{f}_E(\xi)
  := \int_{\Bbb R^n} f(x)\ {\rm e}^{-{\rm i}x\cdot\xi}
  \ \tilde{\rm d}x,
$$
where $\tilde{\rm d}x = (2\pi)^{-n}\ {\rm d}x$.
Then ${\Fcal}_E$ is a bijection and
$$
  f(x) = \int_{\Bbb R^n} \widehat{f}_E(\xi)
  \ {\rm e}^{{\rm i}x\cdot\xi}
  \ {\rm d}\xi,
$$
and this Fourier transform can be uniquely 
extended to
${\Fcal}_E: {\Scal}'(\Bbb R^n)\to{\Scal}'(\Bbb R^n)$.

The main symbol class in the sequel consists of
H\"ormander's $(\rho,\delta)$-symbols of order $m$:
Let $m\in\Bbb R$ and $0\leq\delta<\rho\leq 1$.
For $\xi\in\Bbb R^n$ define $\langle\xi\rangle := (1+\|\xi\|^2)^{1/2}$,
where $\|\xi\|^2 := \sum_{j=1}^n |\xi_j|^2$.
Then
$S^m_{\rho,\delta}(\Bbb R^n\times\Bbb R^n)$
consists of those functions
$\sigma\in C^\infty(\Bbb R^n\times\Bbb R^n)$ for which
\begin{equation}
  \left|\partial_\xi^\alpha \partial_x^\beta \sigma(x,\xi) \right|
        \leq C_{\sigma\alpha\beta m}
                \ \langle\xi\rangle^{m-\rho|\alpha|+\delta|\beta|}
\end{equation}
for every $x\in\Bbb R^n$ and for every $\alpha,\beta\in\Bbb N^n$.

Let $\Bbb T^n = (\Bbb R/2\pi\Bbb Z)^n$ denote the $n$-dimensional torus.
We may identify $\Bbb T^n$ with the hypercube
$[0,2\pi[^n\subset\Bbb R^n$ (or $[-\pi,\pi[^n$).
Functions on $\Bbb T^n$ may be thought as those functions on $\Bbb R^n$
that are $2\pi$-periodic in each of the coordinate directions.
Let ${\Dcal}(\Bbb T^n)$ be the vector space $C^\infty(\Bbb T^n)$
endowed with the usual test function topology,
and let ${\Dcal}'(\Bbb T^n)$ be its dual,
the space of distributions on $\Bbb T^n$.
Inclusion ${\Dcal}(\Bbb T^n)\subset{\Dcal}'(\Bbb T^n)$
is interpreted by
$$
  \phi(\psi) := \int_{\Bbb T^n} \phi(x)\ \psi(x)\ {\rm d}x,
$$
where we identify the measure on torus with the corresponding restriction of
the Euclidean measure on the hypercube.
Let ${\Scal}(\Bbb Z^n)$ denote the space of
{\it rapidly decaying functions $\Bbb Z^n\to\Bbb C$}.
Let ${\Fcal}_T:{\Dcal}(\Bbb T^n)\to{\Scal}(\Bbb Z^n)$
be the {\it toroidal Fourier transform}
(hence the subscript $_T$)
defined by
$$
  ({\Fcal}_T f)(\xi) = \widehat{f}_T(\xi)
  := \int_{\Bbb T^n} f(x)\ {\rm e}^{-{\rm i}x\cdot\xi}
  \ \tilde{\rm d}x,
$$
where $\tilde{\rm d}x = (2\pi)^{-n}\ {\rm d}x$.
Then ${\Fcal}_T$ is a bijection and
$$
  f(x) = \sum_{\xi\in\Bbb Z^n} \widehat{f}_T(\xi)
  \ {\rm e}^{{\rm i}x\cdot\xi}.
$$
This Fourier transform is extended uniquely to
${\Fcal}_T: {\Dcal}'(\Bbb T^n)\to {\Scal}'(\Bbb Z^n)$.
Notice that ${\Scal}'(\Bbb Z^n)$ consists of
those functions $\Bbb Z^n\to\Bbb C$ growing at infinity 
at most polynomially.
Any continuous linear operator
$A:{\Dcal}(\Bbb T^n)\to\Dcal(\Bbb T^n)$ can be presented by a formula
$$
  (Af)(x) = \sum_{\xi\in\Bbb Z^n} \sigma_A(x,\xi)
        \ \widehat{f}(\xi)\ {\rm e}^{{\rm i}x\cdot\xi},
$$
where the unique function $\sigma_A\in C^\infty(\Bbb T^n\times\Bbb Z^n)$
is called the {\it symbol of $A$}:
$$
  \sigma_A(x,\xi) = {\rm e}^{-{\rm i}x\cdot\xi} Ae_\xi(x),
$$
where $e_\xi(x) := {\rm e}^{{\rm i}x\cdot\xi}$.
Notice that when
$\widehat{s_A(x)}_T(\xi) = \sigma_A(x,\xi)$,
the Schwartz kernel $K_A$ of $A$ satisfies
$$
  K_A(x,y) = s_A(x)(x-y)
$$
in the sense of distributions.

Next, in analogy to the classical differential calculus,
we discuss difference calculus,
which is needed when dealing with Fourier series operators.

\section{Difference calculus}

Let $\sigma:\Bbb Z^n\to\Bbb C$.
Let $v_j\in\Bbb N^n$, $(v_j)_j=1$ and $(v_j)_i = 0$
if $i\not=j$.
Let us define the partial difference operator $\triangle_{\xi_j}$ by
$$
  \triangle_{\xi_j} \sigma(\xi) := \sigma(\xi+v_j)-\sigma(\xi),
$$
and define
$$
  \triangle_\xi^{\alpha} := \triangle_{\xi_1}^{\alpha_1}
        \cdots \triangle_{\xi_n}^{\alpha_n}
$$
for $\alpha = (\alpha_j)_{j=1}^n\in\Bbb N^n$.

\begin{lemma} By the binomial theorem,
{\it
$$
  \triangle_\xi^\alpha \sigma(\xi)=\sum_{\beta\leq\alpha} (-1)^{|\alpha-\beta|} {\alpha\choose\beta} \sigma(\xi+\beta).
$$
}
\end{lemma}
By induction, one can show:

\begin{lemma}[Leibnitz formula for differences]
{\it
Let $\phi,\psi:\Bbb Z^n\to\Bbb C$. Then
$$
  \triangle_\xi^\alpha(\phi \psi)(\xi) =
        \sum_{\beta\leq\alpha} {\alpha\choose \beta}
        \left(\triangle_\xi^\beta \phi(\xi)\right)
        \ \triangle_\xi^{\alpha-\beta} \psi(\xi+\beta).
$$
}
\end{lemma}

``Integration by parts'' has the discrete analogy
``summation by parts''
$$
  \sum_{\xi\in\Bbb Z^n} \phi(\xi)\ (\triangle_\xi^\alpha\psi)(\xi)
        = -\sum_{\xi\in\Bbb Z^n} (({\triangle_\xi}^\alpha)^t \phi)(\xi)
                \ \psi(\xi),
$$
where
$(({\triangle_{\xi_j}})^t \phi)(\xi) = \phi(\xi)-\phi(\xi-v_j)$,
provided that the series converge absolutely.

For $\xi\in\Bbb Z^n$ and $\gamma\in\Bbb Z^n$,
let us define
$$
  \xi^{(\gamma)}=\prod_{j=1}^n \xi_j^{(\gamma_j)},
$$
where
$$
  \xi_j^{(\gamma_j)}:=\left\{
  \begin{array}{ll}
        \prod_{i=0}^{\gamma_j-1}\ (\xi_j-i),            & \gamma_j > 0,\\
        1,              &       \gamma_j=0,\\
        \prod_{i=\gamma_j+1}^0 (\xi_j-i)^{-1},                  & \gamma_j<0.
  \end{array}\right.
$$
Then
$$
  \triangle_\xi^\alpha \xi^{(\gamma)}
        = \gamma^{(\alpha)}\ \xi^{(\gamma-\alpha)},
$$
in analogy to
$\partial_\xi^\alpha \xi^\gamma
= \gamma^{(\alpha)}\xi^{\gamma-\alpha}$.

Let us now discuss the discrete version of the Taylor's
theorem. For simplicity, let us consider the one
dimensional case first. 
\begin{theorem}[Discrete Taylor's theorem]
{\it
For a function $\sigma:\Bbb Z\to\Bbb C$,
$$
 \sigma(\xi+\eta)=\sum_{\alpha=0}^{N-1} \frac{1}{\alpha!}
        \ (\triangle_\xi^\alpha \sigma)(\xi)\ \eta^{(\alpha)}
        + R_N(\xi,\eta),
        \quad (\xi,\eta\in\Bbb Z,\ N\in\Bbb N),
$$
where
\begin{eqnarray*}
  |\triangle_\xi^\alpha R_N(\xi,\eta)| & \leq & \left\{
  \begin{array}{ll}
        \frac{1}{N!} \eta^{(N)} \max_{0\leq\kappa<\eta}
                \left| \triangle_\xi^{N+\alpha} \sigma(\xi+\kappa) \right|,
                & \eta\geq N, \\
        0,
                & 0\leq\eta<N, \\
        \frac{1}{N!} \left| \eta^{(N)} \right| \max_{\eta\leq\kappa<0}
                \left| \triangle_\xi^{N+\alpha} \sigma(\xi+\kappa) \right|,
                & \eta<0,
  \end{array}\right.\\
  & \leq & \frac{1}{N!} \left| \eta^{(N)} \right|
                \max_{\kappa\in\{0,\ldots,\eta\}}
                \left| \triangle_\xi^{N+\alpha} \sigma(\xi+\kappa) \right|.
\end{eqnarray*}
}
\end{theorem}

Notice that the estimate above resembles closely the Lagrange form of
the error term in the traditional Taylor theorem:
$$
  \left\{
  \begin{array}{ll}
        f(x+y)=\sum_{j=0}^{N-1} \frac{1}{j!} f^{(j)}(x) y^j + R_N(x,y), & \\
        R_N(x,y) = \frac{1}{N!} f^{(N)}(x+\theta) y^N, &
                \theta\in [\min\{0,y\},\max\{0,y\}].
  \end{array}\right.
$$

\paragraph{\sf Proof.}
First assume that $\eta\geq 0$. Then, by the binomial formula,
\begin{equation}\label{EQ:qqq}
  \sigma(\xi+\eta) = (I+\triangle_\xi)^\eta \sigma(\xi)
        = \sum_{\alpha=0}^\eta {\eta\choose\alpha} \triangle_\xi^\alpha
                \sigma(\xi)
        = \sum_{\alpha=0}^\eta \frac{1}{\alpha!}\ \triangle_\xi^\alpha
                \sigma(\xi)\ \eta^{(\alpha)}.
\end{equation}
Thus $R_N(\xi,\eta)=0$ for $0\leq\eta<N$. Therefore
$$
  \triangle_\eta^\alpha R_N(\xi,\eta)|_{\eta=0}=0,
$$
when $0\leq\alpha<N$.
Now let $\eta$ be an arbitrary integer.
We notice that
$\triangle_\eta^N\eta^{(\alpha)} = \alpha^{(N)}\ \eta^{(\alpha-N)}=0$
for $0\leq\alpha<N$,
so that when we apply $\triangle_\eta^N$, we get
$$
  \triangle_\eta^N \sigma(\xi+\eta)=\triangle_\eta^N R_N(\xi+\eta).
$$
We have hence the Cauchy problem
\begin{eqnarray*}
  \left\{
  \begin{array}{ll}
        \triangle_\eta^N R_N(\xi,\eta)=\triangle_\eta^N \sigma(\xi+\eta),
                & \\ \left.
        \triangle_\eta^\alpha R_N(\xi,\eta) \right|_{\eta=0}=0,
                & 0\leq\alpha\leq N-1.
  \end{array}\right.
\end{eqnarray*}
It is enough to prove the estimate for
$|R_N(\xi,\eta)|$ (i.e. $\alpha=0$).
Let us define $\sigma(\eta):=\eta^{(N)}/N!$. Then
$\triangle_\eta^N \sigma(\eta)=N^{(N)} \eta^{(N-N)}/N! = 1$,
and $\triangle_\xi^\alpha\sigma(\xi)|_{\xi=0}=0$ when $0\leq\alpha<N$,
so by the uniqueness of the solution of the Cauchy problem it has to be
$$
  \left\{
  \begin{array}{ll}
  \sum_{\kappa_N=0}^{-1} \sum_{\kappa_{N-1}=1}^{\kappa_N-1} \cdots
                \sum_{\kappa_1=1}^{\kappa_2-1}\ 1
                = \frac{1}{N!} \eta^{(N)},                      &\eta\geq N, \\
        \sum_{\kappa_N=m}^{-1} \sum_{\kappa_{N-1}=\kappa_N}^{-1} \cdots
                \sum_{\kappa_1=\kappa_2}^{-1}\ 1
                = \frac{1}{N!} \left| \eta^{(N)} \right|,       &\eta<0,
  \end{array}\right.
$$
completing the proof.

\bigskip
Let us now deal with discrete Taylor polynomial -like expansions
for a function $f:\Bbb Z^n\to\Bbb C$.
For $b\in\Bbb N$, let us denote
\begin{equation}\label{sumintegral}
  I_k^b := \sum_{0\leq k < b}\quad {\rm and}\quad
  I_k^{-b} := -\sum_{-b\leq k < 0}.
\end{equation}
It is useful to think of $I_\xi^\theta\cdots$
as of a discrete version of the one-dimensional integral
$\int_0^\theta\cdots {\rm d}\xi$.
In this discrete context,
the difference $\triangle_\xi$ takes the role of
the differential operator ${\rm d}/{\rm d}\xi$.

In the sequel,
we adopt the notational conventions
$$
  I_{k_1}^\theta I_{k_2}^{k_1} \cdots I_{k_\alpha}^{k_{\alpha-1}} 1 =
  \begin{cases}
    1, & {\rm if}\ \alpha = 0,\\
    I_k^\theta 1 , & {\rm if}\ \alpha = 1, \\
    I_{k_1}^\theta I_{k_2}^{k_1} 1, & {\rm if}\ \alpha = 2,
  \end{cases}
$$
and so on.

\begin{lemma}
{\it
If $\theta\in\Bbb Z$ and $\alpha\in\Bbb N$ then
\begin{equation}
  I_{k_1}^\theta I_{k_2}^{k_1} \cdots I_{k_\alpha}^{k_{\alpha-1}} 1
  = \frac{1}{\alpha!}\ \theta^{(\alpha)}.
\end{equation}
}
\end{lemma}

\paragraph{\sf Proof.}
The result follows step-by-step from the observations $k^{(0)} \equiv 1$,
$\triangle_k k^{(i)} = i\ k^{(i-1)}$
and $I_k^b \triangle_k k^{(i)} = b^{(i)}$.

\begin{remark} This is like applying a discrete trivial version
of the fundamental theorem of calculus:
$\int_0^\theta f'(\xi)\ {\rm d}\xi = f(\theta)-f(0)$
for smooth enough $f:\Bbb R\to\Bbb C$
corresponds to
$I_\xi^\theta \triangle_\xi f(\xi) = f(\theta)-f(0)$
for $f:\Bbb Z\to\Bbb C$.
\end{remark}

\begin{corollary}
{\it
If $\theta\in\Bbb Z^n$ and $\alpha\in\Bbb N^n$ then
\begin{equation}\label{counting}
  \prod_{j=1}^n I_{k(j,1)}^\theta I_{k(j,2)}^{k(j,1)}
  \cdots I_{k(j,\alpha_j)}^{k(j,\alpha_j-1)} 1
  = \frac{1}{\alpha!}\ \theta^{(\alpha)},
\end{equation}
where $\prod_{j=1}^n I_j$ means $I_1 I_2\cdots I_n$, where
$I_j := I_{k(j,1)}^{\theta_j} I_{k(j,2)}^{k(j,1)}
\cdots I_{k(j,\alpha_j)}^{k(j,\alpha_j-1)}$.
}
\end{corollary}

We now have
\begin{theorem}
{\it
Let $p:\Bbb Z^n\to\Bbb C$ and
$$
  r_M(\xi,\theta) := p(\xi+\theta) - \sum_{|\alpha| < M}
        \frac{1}{\alpha!}\ \theta^{(\alpha)} \triangle_\xi^\alpha p(\xi).
$$
Then
\begin{equation}\label{remainderestimate}
  \left|\triangle_\xi^\omega r_M(\xi,\theta)\right| \leq
        c_M\ \max_{|\alpha|=M,\ \nu\in Q(\theta)} \left| \theta^{(\alpha)}
        \ \triangle_\xi^{\alpha+\omega} p(\xi+\nu) \right|,
\end{equation}
where $Q(\theta) := \{\nu\in\Bbb Z^n:
        \ \min(0,\theta_j)\leq \nu_j \leq \max(0,\theta_j) \}$.
}
\end{theorem}

\paragraph{\sf Proof.}
For $0\not=\alpha\in\Bbb N^n$,
let us denote $m_\alpha:= \min\{j:\ \alpha_j\not=0\}$.
For $\theta\in\Bbb Z^n$ and $i\in\{1,\ldots,n\}$,
let us define $\nu(\theta,i,k)\in\Bbb Z^n$ by
$$
  \nu(\theta,i,k) := (\theta_1,\ldots,\theta_{i-1},k,0,\ldots,0),
$$
i.e.
$$
  \nu(\theta,i,k)_j =
  \begin{cases}
    \theta_j, & {\rm if}\ 1\leq j < i, \\
    k, & {\rm if}\ j = i, \\
    0, & {\rm if}\ i < j\leq n.
  \end{cases}
$$
We claim that the remainder can be written in the form
\begin{equation}
  r_M(\xi,\theta) = \sum_{|\alpha| = M} r_\alpha(\xi,\theta),
\end{equation}
where for each $\alpha$,
\begin{equation}\label{remainderinduction}
  r_\alpha(\xi,\theta)
  = \prod_{j=1}^n
  I_{k(j,1)}^{\theta_j}
  I_{k(j,2)}^{k(j,1)}\cdots I_{k(j,\alpha_j)}^{k(j,\alpha_j-1)}
    \ \triangle_\xi^\alpha p(\xi+\nu(\theta,m_\alpha,
      k(m_\alpha,\alpha_{m_\alpha})));
\end{equation}
recall (\ref{sumintegral}) and (\ref{counting}).
The proof of (\ref{remainderinduction}) is by induction:
The first remainder term $r_1$ is of the claimed form, since
$$
  r_1(\xi,\theta)
  = p(\xi+\theta)-p(\xi)
  = \sum_{i=1}^n r_{v_i}(\xi,\theta),
$$
where
$$
  r_{v_i}(\xi,\theta) = I_{k}^{\theta_i} \triangle_\xi^{v_i}
  p(\xi+\nu(\theta,i,k));
$$
here $r_{v_i}$ is of the form (\ref{remainderinduction})
for $\alpha = v_i$, $m(\alpha) = i$ and $\alpha_{m_\alpha} = 1$.
So suppose that the claim $(\ref{remainderinduction})$
is true up to order $|\alpha|=M$.
Then
\begin{eqnarray*}
  r_{M+1}(\xi,\theta)
  & = & r_M(\xi,\theta) - \sum_{|\alpha| = M} \frac{1}{\alpha!}
  \ \theta^{(\alpha)}\ \triangle_\xi^\alpha p(\xi) \\
  & = & \sum_{|\alpha| = M} \left( r_\alpha(\xi,\theta) - \frac{1}{\alpha!}
  \ \theta^{(\alpha)}\ \triangle_\xi^\alpha p(\xi) \right) \\
  & = &
    \sum_{|\alpha| = M} \prod_{j=1}^n
  I_{k(j,1)}^{\theta_j}
  I_{k(j,2)}^{k(j,1)}\cdots I_{k(j,\alpha_j)}^{k(j,\alpha_j-1)}
    \\
  & & \triangle_\xi^\alpha \left[ p(\xi+\nu(\theta,m_\alpha,
      k(m_\alpha,\alpha_{m_\alpha}))) - p(\xi) \right],
\end{eqnarray*}
where we used (\ref{remainderinduction}) and
(\ref{counting}) to obtain the last equality.
Combining this to the observation
$$
  p(\xi+\nu(\theta,m_\alpha,
      k)) - p(\xi)
    = \sum_{i=1}^{m_\alpha}
 I_{\ell}^{\nu(\theta,m_\alpha,k)_i}
    \ \triangle_\xi^{v_i} p(\xi+\nu(\theta,i,\ell)),
$$
we get
\begin{eqnarray*}
  r_{M+1}(\xi,\theta)
  & = &
    \sum_{|\alpha| = M} \prod_{j=1}^n
  I_{k(j,1)}^{\theta_j}
  I_{k(j,2)}^{k(j,1)}\cdots I_{k(j,\alpha_j)}^{k(j,\alpha_j-1)}
    \sum_{i=1}^{m_\alpha}
    I_{\ell(i)}^{\nu(\theta,m_\alpha,k(m_\alpha,\alpha_{m_\alpha}))_i} \\
  & & 
  \triangle_\xi^{\alpha+v_i} p(\xi+\nu(\theta,i,\ell(i))) \\
  & = &
  \sum_{|\beta| = M+1}
  \prod_{j=1}^n
  I_{k(j,1)}^{\theta_j}
  I_{k(j,2)}^{k(j,1)}\cdots I_{k(j,\beta_j)}^{k(j,\beta_j-1)} \\
  &  & \triangle_\xi^\beta p(\xi+\nu(\theta,m_\beta,
      k(m_\beta,\beta_{m_\beta})));
\end{eqnarray*}
the last step here is just simple tedious book-keeping.
Thus the induction proof of (\ref{remainderinduction}) is complete.
Finally, let us prove estimate (\ref{remainderestimate}).
By (\ref{remainderinduction}), we obtain
\begin{eqnarray*}
  \left| \triangle_\xi^\omega r_M(\xi,\theta) \right|
  & = & \left| \sum_{|\alpha|=M} \triangle_\xi^\omega r_\alpha(\xi,\theta)
  \right|\\
  & = & \left| \sum_{|\alpha| = M}
  \prod_{j=1}^n
  I_{k(j,1)}^{\theta_j}
  I_{k(j,2)}^{k(j,1)}\cdots I_{k(j,\alpha_j)}^{k(j,\alpha_j-1)} \right. \\
  & & \left.\triangle_\xi^{\alpha+\omega} p(\xi+\nu(\theta,m_\alpha,
      k(m_\alpha,\alpha_{m_\alpha}))) \right| \\
  & \leq & \sum_{|\alpha| = M} \frac{1}{\alpha!} \left|\theta^{(\alpha)}\right|
  \max_{\nu\in Q(\theta)} \left|
     \triangle_\xi^{\alpha+\omega} p(\xi+\nu) \right|,
\end{eqnarray*}
where in the last step we used (\ref{counting}).
The proof is complete.

\begin{remark}
If $n\geq 2$, there are many alternative forms for remainders
$r_\alpha(\xi,\theta)$.
This is due to the fact that there may be many different
shortest discrete step-by-step paths in the space $\Bbb Z^n$
from $\xi$ to $\xi+\theta$.
In the proof above, we chose just one such path,
traveling via the points
$$
  \xi,\quad \xi+\theta_1 v_1,\quad\ldots,\quad
  \xi+\sum_{i=1}^j \theta_i v_i,\quad\ldots,\quad\xi+\theta.
$$
But if $n=1$,
there is just one shortest discrete path
from $\xi\in\Bbb Z$ to $\theta\in\Bbb Z$,
and in that case
$$
  r_M(\xi,\theta) = I_{k_1}^\theta I_{k_2}^{k_1}\cdots I_{k_M}^{k_{M-1}}
  \triangle_\xi^M p(\xi+k_M).
$$
\end{remark}
Notice also that the discrete Taylor theorem presented above
implies the following smooth Taylor result:

\begin{corollary}
{\it
Let $p\in C^\infty(\Bbb R^n)$ and
$$
  r_M(\xi,\theta) := p(\xi+\theta) - \sum_{|\alpha| < M}
        \frac{1}{\alpha!}\ \theta^{\alpha}
        \left( \frac{\partial}{\partial\xi}\right)^\alpha p(\xi).
$$
Then
\begin{equation}
  \left|\partial_\xi^\omega r_M(\xi,\theta)\right| \leq
        c_M\ \max_{|\alpha|=M,\ \nu\in Q_{\Bbb R^n}(\theta)}
        \left| \theta^{\alpha}
        \ \partial_\xi^{\alpha+\omega} p(\xi+\nu) \right|,
\end{equation}
where $Q_{\Bbb R^n}(\theta) := \{\nu\in\Bbb R^n:
        \ \min(0,\theta_j)\leq \nu_j \leq \max(0,\theta_j) \}$.
}
\end{corollary}

\begin{remark}
We see that in the remainder estimates above,
the cubes $Q(\theta)\subset\Bbb Z^n$ and
$Q_{\Bbb R^n}(\theta)\subset\Bbb R^n$
could be replaced by (discrete, resp. continuous) paths
from $0$ to $\theta$;
especially, $Q_{\Bbb R^n}(\theta)$ could be replaced by the straight line
from $0$ to $\theta$.
\end{remark}

\section{Pseudodifferential operators on the torus}

Let $m\in\Bbb R$ and $0\leq\delta<\rho\leq 1$.
Then
$S^m_{\rho,\delta}(\Bbb T^n\times\Bbb Z^n)$
consists of those functions
$\sigma\in C^\infty(\Bbb T^n\times\Bbb Z^n)$ for which
\begin{equation}
  \left|\triangle_\xi^\alpha \partial_x^\beta \sigma(x,\xi) \right|
        \leq C_{\sigma\alpha\beta m}
                \ \langle\xi\rangle^{m-\rho|\alpha|+\delta|\beta|}
\end{equation}
for every $x\in\Bbb T^n$, for every $\alpha,\beta\in\Bbb N^n$.
If $\sigma_A\in S^m_{\rho,\delta}(\Bbb T^n\times\Bbb Z^n)$,
we denote $A\in{\rm Op}S^m_{\rho,\delta}(\Bbb T^n\times\Bbb Z^n)$.
The class $S^m_{\rho,\delta}(\Bbb T^n\times\Bbb T^n\times\Bbb Z^n)$ consists
of the functions $a\in C^\infty(\Bbb T^n\times\Bbb T^n\times\Bbb Z^n)$
such that
\begin{equation}
  \left|\triangle_\xi^\alpha \partial_x^\beta \partial_y^\gamma
        a(x,y,\xi) \right|
        \leq C_{a\alpha\beta\gamma m}
                \ \langle\xi\rangle^{m-\rho|\alpha|+\delta|\beta+\gamma|}
\end{equation}
for every $x,y\in\Bbb T^n$, for every $\alpha,\beta,\gamma\in\Bbb N^n$;
such a function $a$ is called an {\it amplitude} of order $m\in\Bbb R$
of type $(\rho,\delta)$.
Formally we may define
$$
  ({\rm Op}(a)f)(x) := \int_{\Bbb T^n} f(y) \sum_{\xi\in\Bbb Z^n}
        a(x,y,\xi)\ {\rm e}^{{\rm i}(x-y)\cdot\xi}\ \tilde{\rm d}y
$$
for $f\in{\Dcal}(\Bbb T^n)$.

\begin{remark}
On $\Bbb T^n$,
H\"ormander's usual $(\rho,\delta)$-symbol class of order $m\in\Bbb R$
coincides with the class ${\rm Op}S^m_{\rho,\delta}(\Bbb T^n\times\Bbb Z^n)$
\cite{McLean}.
\end{remark}

\begin{lemma}
{\it
Let $a\in S^m_{\rho,\delta}(\Bbb T^n\times\Bbb T^n\times\Bbb Z^n)$,
and define
$$
  a_\alpha(x,y,\xi):=\left( {\rm e}^{{\rm i}(y-x)}-1 \right)^\alpha a(x,y,\xi),
$$
where $\alpha\in\Bbb N^n$. Then
$\triangle_\xi^\alpha a\in
  S^{m-\rho |\alpha|}_{\rho,\delta}(\Bbb T^n\times\Bbb T^n\times\Bbb Z^n)$
and
$Op(a_\alpha)=Op(\triangle_\xi^\alpha a)$.
}
\end{lemma}

\paragraph{\sf Proof.}
Clearly
$\triangle_\xi^\alpha a\in
  S^{m-\rho |\alpha|}_{\rho,\delta}(\Bbb T^n\times\Bbb T^n\times\Bbb Z^n)$.
Now
\begin{eqnarray*}
  ({\rm Op}(a_\alpha)f)(x) & = &
                \int_{\Bbb T^n} f(y) \sum_{\xi\in\Bbb Z^n} a_\alpha(x,y,\xi)
                \ {\rm e}^{{\rm i} (x-y)\cdot\xi}\ \tilde{\rm d}y \\
        & = & \int_{\Bbb T^n} f(y) \left( \sum_{\xi\in\Bbb Z^n}
                \left({\rm e}^{{\rm i}(y-x)}-1\right)^\alpha a(x,y,\xi)
                \ {\rm e}^{{\rm i}(x-y)\cdot\xi}\right)\ \tilde{\rm d}y\\
        & = & \int_{\Bbb T^n} f(y) \left( \sum_{\xi\in\Bbb Z^n} a(x,y,\xi)
                \ ({\triangle}_\xi^\alpha)^t
                {\rm e}^{{\rm i}(x-y)\cdot\xi} \right)\ \tilde{\rm d}y\\
        & = & \int_{\Bbb T^n} f(y) \left( \sum_{\xi\in\Bbb Z^n}
                {\rm e}^{{\rm i}(x-y)\cdot\xi}
                \ \triangle_\xi^\alpha a(x,y,\xi) \right)\ \tilde{\rm d}y.
\end{eqnarray*}
Thus ${\rm Op}(a_\alpha)={\rm Op}(\triangle_\xi^\alpha a)$.

\begin{theorem}
{\it
For every amplitude
$a\in S^m_{\rho,\delta}(\Bbb T^n\times\Bbb T^n\times\Bbb Z^n)$
there exists a unique symbol
$\sigma\in S^m_{\rho,\delta}(\Bbb T^n\times\Bbb Z^n)$
satisfying
${\rm Op}(a)={\rm Op}(\sigma)$,
where
\begin{equation}
  \sigma(x,\xi) \sim \sum_{\alpha\geq 0} \frac{1}{\alpha!}
        \ \triangle_\xi^\alpha
        \ \partial_y^{(\alpha)} a(x,y,\xi)|_{y=x}.
\end{equation}
}
\end{theorem}

\paragraph{Proof:}
essentially the same as in \cite{TurunenVainikko}.

\section{Periodisation}

Let us define the {\it periodisation operator}
$p:{\Scal}(\Bbb R^n)\to{\Dcal}(\Bbb T^n)$ by
$$
  pu (x) := {\Fcal}_T^{-1}(({\Fcal}_E u)|_{\Bbb Z^n})(x).
$$
Let us describe the extension 
of this periodisation
for some nice-enough classes of distributions.
From the proof of the Poisson summation formula
it follows that 
$$
  pu(x) = \sum_{\xi\in\Bbb Z^n} u(x+2\pi\xi).
$$ 
This formula makes sense almost everywhere for
$u\in L^{1}(\Bbb R^n)$.
Indeed, formally
\begin{eqnarray*}
  (pu)(x)
  & = &
  \sum_{k\in\Z^n} {\rm e}^{{\rm i}x\cdot k}
  \int_{\Bbb R^n} {\rm e}^{-{\rm i}y\cdot k}\ u(y)\ {\rm d}y =  
   \int_{\Bbb R^n}
   u(y) \left(\sum_{k\in\Z^n} {\rm e}^{{\rm i}(x-y)\cdot k}\right)
   \ {\rm d}y \\
  & = &
  \int_{\Bbb R^n} u(y)\ \delta_{\Z^n}(2\pi(x-y))\ {\rm d}y
  = \sum_{\xi\in\Z^n} u(x+2\pi \xi).
\end{eqnarray*}
This calculation can be justified in the standard way.
Now $pu\in L^{1}(\Bbb T^n)$ with
$\|pu\|_{L^{1}(\Bbb T^n)}\leq \|u\|_{L^{1}(\Bbb R^n)}$.
Moreover, if $\xi\in\Bbb Z^n$ then
$$
  \widehat{pu}_T(\xi) = \widehat{u}_E(\xi).
$$
Clearly $p:L^{1}(\Bbb R^n)\to L^{1}(\Bbb T^n)$ is a surjection.
We will also use that
$$(pu)(x)=\sum_{\xi\in\Bbb Z^n} {\rm e}^{{\rm i}x\cdot\xi} 
\widehat{pu}_T(\xi)=\sum_{\xi\in\Bbb Z^n} 
{\rm e}^{{\rm i}x\cdot\xi} \widehat{u}_E(\xi).$$

Let us establish basic properties of pseudodifferential
operators with respect to periodisation.
Let us call symbol
$a(x,\xi)$ periodic if the function $x\mapsto a(x,\xi)$ is
$2\pi$-periodic. We will use tildes to denote
corresponding restricted operators.

\begin{prop}
Let $a\in S^m_{\rho,\delta}(\Rn\times\Rn)$ be a periodic symbol.
Let $\tilde{a}=a|_{\Tn\times\Z^n}$.
Then $p\circ a(X,D)=\tilde{a}(X,D)\circ p.$
\label{p:p1}
\end{prop}

Notice that $a$ does not have to be a symbol, as
the same property holds when we define $a(X,D)$ 
in the usual sense, by even quite irregular amplitude $a(x,\xi)$.

\paragraph{\sf Proof.}
Notice that $\tilde{a}\in S^m_{\rho,\delta}(\Tn\times\Bbb Z^n).$
Let $f\in L^{1}(\Rn)$. Then we have
\begin{eqnarray*}
  p(a(X,D)f)(x) & = & 
  \sum_{k\in\Z^n} a(x+2\pi k,D) f(x+2\pi k) \\
  & = &
  \sum_{k\in\Z^n} \int_{\Rn} {\rm e}^{{\rm i}(x+2\pi k)\cdot\xi}
  \ a(x+2\pi k,\xi)\ \widehat{f}_E(\xi) \ {\rm d}\xi \\
  & = &
  \int_{\Rn} \left(\sum_{k\in\Z^n} {\rm e}^{{\rm i}2\pi k\cdot\xi}
  \right) {\rm e}^{{\rm i}x\cdot\xi}
  \ a(x,\xi)\ \widehat{f}_E(\xi) \ {\rm d}\xi \\
  & = &
  \int_{\Rn}
  {\rm e}^{{\rm i}x\cdot\xi}\ a(x,\xi)\ \widehat{f}_E(\xi)
  \ \delta_{\Z^n}(\xi) \ {\rm d}\xi \\
  & = &
  \sum_{\xi\in\Z^n} {\rm e}^{{\rm i}x\cdot\xi}
  \ a(x,\xi)\ \widehat{f}_E(\xi)  \\
  & = &
  \sum_{\xi\in\Z^n} {\rm e}^{{\rm i}x\cdot\xi}
  \ a(x,\xi)\ \widehat{pf}_T(\xi) \\
  & = &
  \tilde{a}(X,D)(pf)(x);
\end{eqnarray*}
these calculations are justified in the
sense of distributions. The proof is complete.

Since we will not always work with periodic symbols it
may be convenient to periodize them. If $a(X,D)$ is a
pseudodifferential operator with symbol $a(x,\xi)$, 
by $(pa)(X,D)$ we will denote a pseudodifferential
operator with symbol 
$(pa)(x,\xi)=\sum_{k\in\Z^n} a(x+2\pi k,\xi)$. 
This makes sense if, for example,  
$a$ in integrable in $x$.

\begin{prop}
Let $a\in S^m_{\rho,\delta}(\Rn\times\Rn)$ satisfy
$a(x,\xi)=0$ for all $x\in\Bbb R^n\setminus [-\pi,\pi]^n$.
Then we have
$$
a(X,D)f=(pa)(X,D)f+Rf,
$$
for all $f$ supported in $[-\pi,\pi]^n$. Here 
$R:{\Scal}^\prime(\Rn)\to {\Scal}(\Rn)$ is a smoothing
pseudodifferential operator.
\label{p:p2}
\end{prop}

\paragraph{\sf Proof.}
By our definition we can write
$$
  (pa)(X,D)f(x)=\sum_{k\in\Z^n}\int_{\Bbb R^n} {\rm e}^{{\rm i}x\cdot\xi}
  \ a(x+2\pi k,\xi)\ \widehat{f}_E(\xi) \ {\rm d}\xi,
$$
and let $Rf=a(X,D)f-(pa)(X,D)f.$
The assumption on the support of $a$ implies
that for every $x$ there is only one $k\in\Z^n$ for which
$a(x+2\pi k,\xi)\not=0$, so the sum consists of only one
term. It follows that $Rf(x)=0$ for $x\in [-\pi,\pi]^n$. 
Let now $x\in\Bbb R^n\setminus [-\pi,\pi]^n$. Since
$$
  Rf(x) = -\sum_{k\in\Z^n, k\not=0} \int_{\Bbb R^n} \int_{\Bbb R^n}
  {\rm e}^{{\rm i}(x-y)\cdot\xi}
  \ a(x+2\pi k,\xi)\ f(y)\ {\rm d}y\ {\rm d}\xi
$$
is just a single term and $|x-y|>0$, we can integrate by
parts with respect to $\xi$ any number of times. This
implies that $R\in\Psi^{-\infty}$ and that $Rf$ decays
at infinity faster than any power. The proof is complete
since the same argument can be applied to the derivatives of
$Rf$. 

\begin{remark}
Note that if $f$ is compactly supported, but not necessarily 
in the cube $[-\pi,\pi]^n$, sums in the proof may consist of
finite number of terms. This means that on ${\Ecal}^\prime(\Rn)$,
modulo a smoothing operator,
we can write $a(X,D)$ as a finite sum  of operators with 
periodic symbols. Moreover, the same argument applies if
$a(x,\xi)$ is compactly supported in $x$, but not necessarily
in $[-\pi,\pi]^n$.
\end{remark}

This proposition allows us to extend formula of Proposition
\ref{p:p1} to perturbations of periodic symbols. We will use
it when $a(x,D)$ is a sum of a constant coefficient operator
and an operator with symbol having compact $x$-support.

\begin{corollary}
Let $a(X,D)$ be an operator with symbol
$$a(x,\xi)=a_1(x,\xi)+a_0(x,\xi),$$
where $a_1\in S^m_{\rho,\delta}(\Rn\times\Rn)$ is periodic in $x$ and 
$a_0\in S^m_{\rho,\delta}(\Rn\times\Rn)$ has compact $x$-support.
Then there is
a symbol $\tilde{b}\in S^m_{\rho,\delta}(\Tn\times\Zn)$ such that
$$
  p(a(X,D)f)=\tilde{b}(X,D)(pf)+{p}(Rf),
  \; f\in {\Ecal}^\prime(\Rn),
$$
where $R:{\Scal}^\prime(\Rn)\to{\Scal}(\Rn).$
In particular, if $\supp(a_0(\cdot,\xi)), 
\supp(f)\subset [-\pi,\pi]^n$, we can take
$\tilde{b}(X,D)=\widetilde{a_1}(X,D)+\widetilde{pa_0}(X,D).$
\label{p:p3}
\end{corollary}  
We assumed that symbols are smooth, but the
requirement of the smoothness of $a_1(x,\xi)$ is not
necessary similar to Proposition \ref{p:p1}.
\paragraph{\sf Proof.} 
By Proposition \ref{p:p2},
$a(X,D)=a_1(X,D)+(pa_0)(X,D)+R.$ Since operator
$b(X,D)=a_1(X,D)+(pa_0)(X,D)$ has periodic symbol,
by Proposition \ref{p:p1} we have
$p\circ b(X,D)=\tilde{b}(X,D)\circ p=\widetilde{a_1}(X,D)\circ p+
\widetilde{pa_0}(X,D)\circ p.$
Since $R:{\Scal}^\prime(\Rn)\to{\Scal}(\Rn)$, we also have
$p\circ R:{\Scal}^\prime(\Rn)\to{\Dcal}(\Bbb T^n).$ The proof
is complete.

\section{Conditions for $L^{2}$-boundedness}

Next we study conditions on a toroidal symbol $\sigma_A$
that guarantee $L^{2}$-boundedness for the corresponding operator
$A:{\Dcal}(\Bbb T^n)\to{\Dcal}(\Bbb T^n)$.
Notice that $x\mapsto\sigma_A(x,\xi)\in C^\infty(\Bbb T^n)$
for every $\xi\in\Bbb Z^n$.

\begin{prop}
{\it
If
$$
  \left| \partial_x^\beta \sigma_A(x,\xi) \right| \leq
  C
$$
when $|\beta| \leq n/2+1$
then $A\in{\Lcal}(L^{2}(\Bbb T^n))$.
}
\end{prop}

\paragraph{\sf Proof.}
Now
\begin{eqnarray*}
  Af(x)
  & = & \sum_{\xi\in\Bbb Z^n} \sigma_A(x,\xi)\ \widehat{f}_T(\xi)
  \ {\rm e}^{{\rm i}x\cdot\xi} \\
  & = & \sum_{\xi,\eta\in\Bbb Z^n} \widehat{\sigma_A}_T(\eta,\xi)
  \ \widehat{f}_T(\xi)
  \ {\rm e}^{{\rm i}x\cdot(\xi+\eta)}\\
  & = & \sum_{\omega\in\Bbb Z^n} {\rm e}^{{\rm i}x\cdot\omega}
  \sum_{\xi\in\Bbb Z^n} \widehat{\sigma_A}_T(\omega-\xi,\xi)
  \ \widehat{f}_T(\xi).
\end{eqnarray*}
Here
$\left|\widehat{\sigma_A}_T(\eta,\xi)\right| \leq C\ \langle\eta\rangle^{-k}$,
so that
\begin{eqnarray*}
  \| A f \|_{L^{2}(\Bbb T^n)}^2
  & = & \int_{\Bbb T^n} |Af(x)|^2\ \tilde{\rm d}x \\
  & = & \sum_{\omega\in\Bbb Z^n} |\widehat{Af}_T(\omega)|^2 \\
  & = & \sum_{\omega\in\Bbb Z^n} \left|
    \sum_{\xi\in\Bbb Z^n} \widehat{\sigma_A}_T(\omega-\xi,\xi)
    \ \widehat{f}_T(\xi) \right|^2 \\
  &\leq &
  \left( \sup_{\omega\in\Bbb Z^n} \sum_{\xi\in\Bbb Z^n}
    \left|\widehat{\sigma_A}_T(\omega-\xi,\xi)\right|^2
    \right) \\
  & & \left( \sup_{\xi\in\Bbb Z^n} \sum_{\omega\in\Bbb Z^n}
    \left|\widehat{\sigma_A}_T(\omega-\xi,\xi)\right|^2
    \right)
    \sum_{\xi\in\Bbb Z^n} \left|\widehat{f}_T(\xi) \right|^2 \\
  & \leq & C'\ \|f\|_{L^{2}(\Bbb T^n)}^2.
\end{eqnarray*}

Note that no difference conditions for the $\xi$-variable 
were needed.
In fact, this is related to the following more
general result.

\begin{theorem}\label{th:fiol2}
{\it
Let $$Tu(x)=\sum_{k\in\Zn} {\rm e}^{{\rm i}\phi(x,k)}\ a(x,k)\ \hat{u}(k).$$
Assume that for every $\alpha$ for which $|\alpha|\leq 2n+1$,
\begin{equation}\label{l2a1}
  \left| \partial_x^\alpha a(x,k) \right| \leq
  C,\quad\quad\quad
  \left| \partial_x^\alpha\Delta_k \phi(x,k) \right| \leq
  C.
\end{equation}
Assume also that for every $x\in\Tn$ and for every $k,l\in\Zn$,
\begin{equation}\label{l2a2}
  \left| \nabla_x\phi(x,k)-\nabla_x\phi(x,l) \right| \geq
  C|k-l|.
\end{equation}
Then $T\in{\mathcal L}(L^{2}(\Bbb T^n))$.
}
\end{theorem}

Note that condition \eqref{l2a2} is a discrete version of
the usual local graph condition for Fourier integral operators,
necessary for the local $L^2$-boundedness. We also note that
these conditions roughly correspond to $C^1$ properties of the phase
in $\xi$.
Finally, if $\phi$ and $a$ are not $2\pi$--periodic in $x$,
operator $T$ is bounded from $L^2(\Tn)$ to $L^2_{loc}(\Rn)$.
Theorem \ref{th:fiol2} is the discrete version of the global
boundedness theorem in \cite{RSCPDE}.

\section{Extending symbols}\label{extending}

It is often useful to extend toroidal symbols
from $\Bbb T^n\times\Bbb Z^n$
to $\Bbb T^n\times\Bbb R^n$.
This can be done with a suitable convolution
that respects the symbol inequalities.
The idea of the following lemma goes probably back to Y. Meyer.

\begin{lemma}\label{l:aux}
{\it
There exist
$\theta,\phi_\alpha\in{\Scal}(\Bbb R^n)$ such that
$\widehat{\theta}_E|_{\Bbb Z^n}(\xi) = \delta_{0,\xi}$ and
$\partial_\xi^\alpha\widehat{\theta}_E(\xi)
= (\triangle_\xi^\alpha)^t \phi_\alpha(\xi)$
for every multi-index $\alpha$.
}
\end{lemma}

\paragraph{\sf Proof.}
Let us first consider the case $n=1$.
Let $\theta=\theta_1\in C^\infty(\Bbb R^{1})$ such that
$$
  {\rm supp}(\theta_1) \subset ]-2\pi,2\pi[,\quad
  \theta_1(-x) = \theta_1(x),\quad
  \theta_1(\pi-y) + \theta_1(\pi+y) = 1
$$
for $x\in\Bbb R$ and
for $0\leq y\leq\pi$.
These assumptions for $\theta$ are enough for us,
and of course the choice is not unique.
In any case, $\widehat{\theta_1}_E \in {\Scal}(\Bbb R^{1})$.
If $\xi\in\Bbb Z^n$ then
\begin{eqnarray*}
  \widehat{\theta_1}_E(\xi)
  & = &
  \int_{\Bbb R^{1}} \theta_1(x)\ {\rm e}^{-{\rm i}x\cdot\xi}
  \ \tilde{\rm d}x \\
  & = &
  \int_0^{2\pi} \left(\theta_1(x-2\pi) + \theta_1(x) \right)
  \ {\rm e}^{-{\rm i}x\cdot\xi}\ \tilde{\rm d}x \\
  & = &
  \int_0^{2\pi}
  \ {\rm e}^{-{\rm i}x\cdot\xi}\ \tilde{\rm d}x \\
  & = &
  \delta_{0,\xi}.
\end{eqnarray*}
If desired $\phi_\alpha\in{\Scal}(\Bbb R^n)$ exists,
it must satisfy
$$
  \int_{\Bbb R^{1}} {\rm e}^{{\rm i}x\cdot\xi}\ \phi^{(\alpha)}(\xi)
  \ {\rm d}\xi =
  \int_{\Bbb R^{1}} {\rm e}^{{\rm i}x\cdot\xi}
  \ (\triangle_\xi^\alpha)^t \phi_\alpha(\xi)
  \ {\rm d}\xi
$$
because
${\Fcal}_E:{\Scal}(\Bbb R^{1})\to{\Scal}(\Bbb R^{1})$
is bijective.
Integration by parts yields
$$
  (-{\rm i}x)^\alpha \theta(x) = (1-{\rm e}^{{\rm i}x})^\alpha
  ({\Fcal}_E^{-1}\phi_\alpha)(x).
$$
Thus
$$
  ({\Fcal}_E^{-1}\phi_\alpha)(x) =
  \left( \frac{-{\rm i}x}{1-{\rm e}^{{\rm i}x}} \right)^{\alpha}
  \theta(x).
$$
The general $n$-dimensional case is reduced to the $1$-dimensional case,
since
$\theta=(x\mapsto\theta_1(x_1)\theta_1(x_2)\cdots \theta_1(x_n))
\in{\Scal}(\Bbb R^n)$
has desired properties.
The following two results can be easily obtained from the
discrete Taylor's theorem.

\begin{lemma}
{\it
Let $a:\Bbb T^n\times\Bbb R^n\to\Bbb C$
belong to $a\in S^m_{\rho,\delta}(\Bbb R^n\times\Bbb R^n)$.
Then the restriction
$\sigma
= a|_{\Bbb T^n\times\Bbb Z^n}\in S^m_{\rho,\delta}(\Bbb T^n\times \Bbb Z^n)$.
}
\end{lemma}

\begin{prop}\label{p:uni}
{\it
Let $a,b:\Bbb T^n\times\Bbb R^n\to\Bbb C$ such that
$a,b\in S^m_{\rho,\delta}(\Bbb R^n\times\Bbb R^n)$
and
$a|_{\Bbb T^n\times\Bbb Z^n} = b|_{\Bbb T^n\times\Bbb Z^n}$.
Then $a-b$ is smoothing,
$a-b\in S^{-\infty}(\Bbb R^n\times\Bbb R^n)$.
}
\end{prop}

The main theorem of this paragraph is that we can extend 
toroidal symbols in a unique smooth way.

\begin{theorem}
{\it
Let $\sigma\in S^m_{\rho,\delta}(\Bbb T^n\times\Bbb Z^n)$.
Then there exists $a\in S^m_{\rho,\delta}(\Bbb R^n\times\Bbb R^n)$
such that
$\sigma = a|_{\Bbb T^n\times\Bbb Z^n}$;
this extended symbol is unique up to smoothing.
}
\end{theorem}

\paragraph{\sf Proof.}
Uniqueness up to smoothing follows from Proposition 
\ref{p:uni},
so the existence is the main issue here.
Let $\theta\in{\Scal}(\Bbb R^n)$ be as in Lemma 
\ref{l:aux}.
Define $a:\Bbb R^n\times\Bbb R^n\to\Bbb C$ by
$$
  a(x,\xi) := \sum_{\eta\in\Bbb Z^n} \widehat{\theta}_E(\xi-\eta)
  \ \sigma(x,\eta).
$$
It is easy to see that $\sigma=a|_{\Bbb T^n\times\Bbb Z^n}$.
Furthermore,
\begin{eqnarray*}
  \left|\partial_\xi^\alpha\partial_x^\beta a(x,\xi)\right|
  & = &
  \left| \sum_{\eta\in\Bbb Z^n}
    \partial_\xi^\alpha\widehat{\theta}_E(\xi-\eta)
    \ \partial_x^\beta \sigma(x,\eta)
  \right| \\
  & = &
  \left| \sum_{\eta\in\Bbb Z^n}
    (\triangle_\xi^\alpha)^t \phi_\alpha(\xi-\eta)
    \ \partial_x^\beta \sigma(x,\eta)
  \right| \\
  & = &
  \left| \sum_{\eta\in\Bbb Z^n}
    \phi_\alpha(\xi-\eta)
    \ {\triangle}_\eta^\alpha\partial_x^\beta \sigma(x,\eta)
    \ (-1)^{|\alpha|} \right| \\
  & \leq &
    \sum_{\eta\in\Bbb Z^n}
    |\phi_\alpha(\xi-\eta)|
    \ C_{\alpha\beta m}\ \langle\eta\rangle^{m-\rho|\alpha|+\delta|\beta|} \\
  & \leq &
    C_{\alpha\beta m}'\ \langle\xi\rangle^{m-\rho|\alpha|+\delta|\beta|}
    \sum_{\eta\in\Bbb Z^n} |\phi_\alpha(\eta)|
    \ \langle\eta\rangle^{\left|m-\rho|\alpha|+\delta|\beta|\right|} \\
  & \leq &
    C_{\alpha\beta m}''\ \langle\xi\rangle^{m-\rho|\alpha|+\delta|\beta|}.
\end{eqnarray*}
Thus $a\in S^m_{\rho,\delta}(\Bbb R^n\times\Bbb R^n)$.

\section{Fourier series operator calculus}

In this section we will describe composition formulae
of Fourier series operators with pseudo-differential operators.
They are similar to the global composition formulae in
\cite{RuzhanskySugimoto1} and \cite{RuzhanskySugimoto}
in $\Rn$. However, the situation on the torus is technically
much simpler since it does not require the global in space
analysis of the corresponding remainders.

\begin{theorem}[composition $TP$]
{\it
Let $T:{\Dcal}(\Bbb T^n)\to{\Dcal}'(\Bbb T^n)$ be defined by
$$
  Tu(x) := \sum_{\xi\in\Bbb Z^n} \int_{\Bbb T^n}
  {\rm e}^{{\rm i}(\phi(x,\xi)-y\cdot\xi)}
  \ a(x,y,\xi)\ u(y)
  \ \tilde{\rm d}y,
$$
where
the amplitude
$a\in C^\infty(\Bbb T^n\times\Bbb T^n\times\Bbb Z^n)$
satisfies
$$
  \left| \partial_x^\alpha \partial_y^\beta a(x,y,\xi) \right|
  \leq C_{\alpha\beta m}\ \langle\xi\rangle^m
$$
for every $x,y\in\Bbb T^n$, $\xi\in\Bbb Z^n$ and
$\alpha,\beta\in\Bbb N^n$;
no restrictions for $\phi$ here.
Let $p\in S^t(\Bbb T^n\times\Bbb Z^n)$.
Then
$$
  TPu(x) = \sum_{\xi\in\Bbb Z^n} \int_{\Bbb T^n}
  {\rm e}^{{\rm i}(\phi(x,\xi)-z\cdot\xi)}
  \ c(x,z,\xi)\ u(z)
  \ \tilde{\rm d}z,
$$
where
$$
  c(x,z,\xi) = \sum_{\eta\in\Bbb Z^n} \int_{\Bbb T^n}
  {\rm e}^{{\rm i}(y-z)\cdot(\eta-\xi)}
  \ a(x,y,\xi)\ p(y,\eta)
  \ \tilde{\rm d}y
$$
satisfying
$$
  \left| \partial_x^\alpha \partial_z^\beta c(x,z,\xi) \right|
  \leq C_{\alpha\beta m t}\ \langle\xi\rangle^{m+t}
$$
for every $x,z\in\Bbb T^n$, $\xi\in\Bbb Z^n$ and
$\alpha,\beta\in\Bbb N^n$.
Moreover,
$$
  c(x,z,\xi) \sim \sum_{\alpha\in\geq 0}
  \frac{1}{\alpha!}\ \partial_y^{(\alpha)}
  \left[ a(x,y,\xi)\ \triangle_\xi^\alpha p(y,\xi) \right] |_{y=z}.
$$
}
\end{theorem}

Composition in the other direction is given by the following
theorem.
\begin{theorem}[composition $PT$]
{\it
Let $T:{\Dcal}(\Bbb T^n)\to{\Dcal}'(\Bbb T^n)$ such that
$$
  Tu(x) := \sum_{\xi\in\Bbb Z^n} \int_{\Bbb T^n}
  {\rm e}^{{\rm i}(\phi(x,\xi)-y\cdot\xi)}
  \ a(x,y,\xi)\ u(y)
  \ \tilde{\rm d}y,
$$
where $a\in C^\infty(\Bbb T^n\times\Bbb T^n\times\Bbb Z^n)$ satisfying
$$
  \left|\partial_x^\alpha\partial_y^\beta a(x,y,\xi) \right|
  \leq C_{\alpha\beta m}\ \langle\xi\rangle^m
$$
for every $x,y\in\Bbb T^n$, $\xi\in\Bbb Z^n$ and
$\alpha,\beta\in\Bbb N^n$;
we assume that
$\phi\in C^\infty(\Bbb T^n\times\Bbb Z^n)$ satisfies
$$
  C^{-1}\ \langle\xi\rangle
  \leq \langle\nabla_x\phi(x,\xi)\rangle
  \leq C\ \langle\xi\rangle
$$
for some $C$ for every $x\in\Bbb T^n$, $\xi\in\Bbb Z^n$,
and that
$$
  \left|\partial_x^\alpha\phi(x,\xi)\right| \leq C_\alpha\ \langle\xi\rangle,
  \quad
  \left|\partial_x^\alpha \triangle_\xi^\beta \phi(x,\xi) \right|
  \leq C_{\alpha\beta}
$$
for every $x\in\Bbb T^n$, $\xi\in\Bbb Z^n$ and
$\alpha,\beta\in\Bbb N^n\setminus\{0\}$.
Let $p\in S^t(\Bbb T^n\times\Bbb Z^n)$.
Then
$$
  p(x,D)Tu(x) = \sum_{\xi\in\Bbb Z^n} \int_{\Bbb R^n}
  {\rm e}^{{\rm i}(\phi(x,\xi)-z\cdot\xi)}
  \ c(x,z,\xi)\ u(z)
  \ \tilde{\rm d}z,
$$
where
$$
  \left| \partial_x^\alpha\partial_z^\beta c(x,z,\xi) \right|
  \leq C_{\alpha\beta}\ \langle\xi\rangle^{m+t}
$$
for every $x,z\in\Bbb T^n$, $\xi\in\Bbb Z^n$ and
$\alpha,\beta\in\Bbb N^n$.
Moreover,
$$
  c(x,z,\xi) \sim \sum_{\alpha\geq 0}
  \frac{{\rm i}^{-|\alpha|}}{\alpha!}
    \ \partial_\eta^\alpha p(x,\eta)|_{\eta=\nabla_x\phi(x,\xi)}
    \ \partial_y^\alpha \left[
      {\rm e}^{{\rm i}\Psi(x,y,\xi)} a(y,z,\xi) \right] |_{y=x},
$$
(here we use a smooth extension for the symbol $p(x,\eta)$)
where
\be \label{eq:Psi}
  \Psi(x,y,\xi) := \phi(y,\xi) - \phi(x,\xi) + (x-y)\cdot\nabla_x\phi(x,\xi),
\ee
when $x\approx y$.
}
\end{theorem}


\section{Applications to hyperbolic equations}

Let $a(x,D)\in\Psi^m(\Bbb R^n)$
(with some properties to be specified).
If $u$ depends on $x$ and $t$, we write
\begin{eqnarray*}
  a(x,D) u(x,t) & = & \int_{\Bbb R^n} a(x,\xi)
  \ \widehat{u}_E(\xi,t)\ {\rm e}^{{\rm i}x\cdot\xi}\ {\rm d}\xi \\
  & = & \int_{\Bbb R^n} \int_{\Bbb R^n} {\rm e}^{{\rm i}(x-y)\cdot\xi}
  \ a(x,\xi)\ u(y,t)\ \tilde{\rm d}y\ {\rm d}\xi.
\end{eqnarray*}
Let $u(\cdot,t)\in L^{1}(\Bbb R^n)$ ($0<t<t_0$) be a solution
to the hyperbolic problem
\be \label{he1}
  \begin{cases}
  {\rm i} \frac{\partial}{\partial t}u(x,t) = a(x,D) u(x,t), \\
  u(x,0) = f(x),
  \end{cases}
\ee
where $f\in L^{1}(\Bbb R^n)$ is compactly supported.

Assume now that $a(X,D)=a_1(X,D)+a_0(X,D)$ where $a_1(x,\xi)$
is periodic and $a_0(x,\xi)$ is compactly supported in $x$
(assume even that $\supp(a_0(\cdot,\xi))\subset [-\pi,\pi]^n$).
Typically, we will want to have $a_1(x,\xi)=a_1(\xi)$ a 
constant coefficient operator, not necessarily smooth in $\xi$.
Let us also assume that $\supp(f) \subset [-\pi,\pi]^n$.

We will now describe a way to periodise problem \eqref{he1}.
According to Proposition \ref{p:p2}, we can replace \eqref{he1}
by 
\be \label{he2}
  \begin{cases}
  {\rm i} \frac{\partial}{\partial t}u(x,t) = 
  (a_1(x,D)+(p a_0)(X,D)) u(x,t)+Ru(x,t), \\
  u(x,0) = f(x),
  \end{cases}
\ee
where the symbol $a_1+pa_0$ is periodic and $R$ is a smoothing
operator. To study singularities of \eqref{he1},
it is sufficient to analyse the Cauchy problem
\be \label{he3}
  \begin{cases}
  {\rm i} \frac{\partial}{\partial t}v(x,t) = 
  (a_1(x,D)+(p a_0)(X,D)) v(x,t), \\
  v(x,0) = f(x)
  \end{cases}
\ee
since by Duhamel's formula ${\rm WF}(u-v)=\emptyset$.
This problem can be transfered to the torus. 
Let $w(x,t)=pv(\cdot,t)(x)$. In view of Proposition \ref{p:p1}
it will solve Cauchy problem
\be \label{he4}
  \begin{cases}
  {\rm i} \frac{\partial}{\partial t}w(x,t) = 
  (\widetilde{a_1}(x,D)+\widetilde{p a_0}(X,D)) w(x,t), \\
  w(x,0) = pf(x).
  \end{cases}
\ee
Calculus constructed in previous sections provides the solution
in the form
$$
  w(x,t)=\sum_{k\in\Z^n} {\rm e}^{{\rm i}\phi(t,x,k)} c(t,x,k)
  \widehat{f}_E(k).
$$
Here we note that $\widehat{pf}_T(k)=\widehat{f}_E(k)$. Also, if 
the symbol $a_1(x,\xi)=a_1(\xi)$ has constant coefficients
and $a_0$ is of order zero, we
have $\phi(t,x,k)=x\cdot k+t a_1(k).$ In particular,
$\nabla_x\phi(x,k)=k$, so composition formulas for $b(x,D)w$
in previous sections can be applied. Details of this 
analysis and investigation of the corresponding properties will
appear in \cite{RT}.

\par


\begin{thebibliography}{99}

\bibitem{Agranovich1}
M. S. Agranovich,
{\it Spectral properties of elliptic pseudodifferential operators
on a closed curve.}
(Russian) Funktsional. Anal. i Prilozhen.
{\bf 13} (1979), no. 4, 54--56.

\bibitem{Agranovich2}
M. S. Agranovich,
{\it Elliptic pseudodifferential operators on a closed curve.}
(Russian) Trudy Moskov. Mat. Obshch.
{\bf 47} (1984), 22--67, 246. 

\bibitem{Amosov}
B. A. Amosov,
{\it On the theory of pseudodifferential operators on the circle.}
(Russian) Uspekhi Mat. Nauk {\bf 43} (1988), no. 3(261), 169--170;
translation in Russian Math. Surveys {\bf 43} (1988), no. 3, 197--198.

\bibitem{Boulkhemair}
A. Boulkhemair, 
{\it $L\sp 2$ estimates for pseudodifferential operators.} 
Ann. Scuola Norm. Sup. Pisa Cl. Sci. (4) {\bf 22} (1995), 
no. 1, 155--183. 

\bibitem{Elschner}
J. Elschner,
{\it Singular ordinary differential operators
and pseudodifferential equations.}
Lecture Notes in Mathematics, 1128. Springer-Verlag, Berlin, 1985.

\bibitem{GarelloMorando}
G. Garello, A. Morando, 
{\it $L\sp p$-boundedness for pseudodifferential operators 
with non-smooth symbols and applications.} 
Boll. Unione Mat. Ital. Sez. B Artic. Ric. Mat. (8) 
{\bf 8} (2005), no. 2, 461--503. 

\bibitem{Hormander}
L. H\"ormander,
{\it The Analysis of Linear Partial Differential Operators IV.}
Springer-Verlag, 1985.

\bibitem{Kumanogonagase}
H. Kumano-go, M. Nagase, 
{\it Pseudo-differential operators with non-regular symbols 
and applications.} Funkcial. Ekvac. {\bf 21} (1978), no. 2, 151--192. 

\bibitem{McLean}
W. McLean,
{\it Local and global description of periodic pseudodifferential operators.}
Math. Nachr. {\bf 150} (1991), 151--161.

\bibitem{ProssdorfSchneider}
S. Pr\"ossdorf, R. Schneider,
{\it Spline approximation methods for multidimensional
periodic pseudodifferential equations.}
Integral Equations Operator Theory {\bf 15} (1992), no. 4, 626--672. 

\bibitem{RuzhanskySugimoto1}
M. Ruzhansky, M. Sugimoto,
{\it Global calculus of Fourier integral operators, 
weighted estimates, and applications to global analysis of 
hyperbolic equations}, 
in Advances in pseudo-differential operators, 65--78, 
Oper. Theory Adv. Appl., 164, Birkh\"auser, 2006.

\bibitem{RSCPDE}
M. Ruzhansky, M. Sugimoto,
{\it Global $L^{2}$ boundedness theorems for a class of 
Fourier integral operators,} to appear in Comm. Partial
Differential Equations.

\bibitem{RuzhanskySugimoto}
M. Ruzhansky, M. Sugimoto,
{\it Weighted $L^{2}$ estimates for a class of Fourier 
integral operators,} preprint.

\bibitem{RT}
M. Ruzhansky, V. Turunen, 
{\it Fourier integral operators on the torus},
in preparation.

\bibitem{SaranenWendland}
J. Saranen, W. L. Wendland,
{\it The Fourier series representation of pseudodifferential operators
on closed curves.}
Complex Variables Theory Appl. {\bf 8} (1987), no. 1-2, 55--64.

\bibitem{Sugimoto}
M. Sugimoto, 
{\it Pseudo-differential operators on Besov spaces.}
Tsukuba J. Math. {\bf 12} (1988), no. 1, 43--63.

\bibitem{Turunen}
V. Turunen,
{\it Commutator characterization of periodic pseudodifferential operators.}
Z. Anal. Anw. {\bf 19} (2000), 95--108. 

\bibitem{TurunenVainikko}
V. Turunen, G. Vainikko,
{\it On symbol analysis of periodic pseudodifferential operators.}
Z. Anal. Anw. {\bf 17} (1998), 9--22. 

\bibitem{VainikkoLifanov1}
G. M. Vainikko, I. K. Lifanov,
{\it Generalization and use of the theory of pseudodifferential operators
in the modeling of some problems in mechanics.} (Russian) 
Dokl. Akad. Nauk {\bf 373} (2000), no. 2, 157--160.

\bibitem{VainikkoLifanov2}
G. M. Vainikko, I. K. Lifanov,
{\it The modeling of problems in aerodynamics and wave diffraction
and the extension of Cauchy-type integral operators on closed and open curves.}
(Russian) 
Differ. Uravn. {\bf 36} (2000), no. 9, 1184--1195, 1293;
translation in Differ. Equ. {\bf 36} (2000), no. 9, 1310--1322.

\end{thebibliography}
\end{document}